\documentclass{gtmon_a}
\pdfoutput=1

\proceedingstitle{The interaction of finite-type and Gromov--Witten
invariants (BIRS 2003)}
\conferencestart{15 November 2003}
\conferenceend{20 November 2003}
\conferencename{The interaction of finite-type and Gromov--Witten
invariants}
\conferencelocation{Banff International Research Station, Banff, Alberta,
Canada}

\editor{David Auckly}
\givenname{David}
\surname{Auckly}

\editor{Jim Bryan}
\givenname{Jim}
\surname{Bryan}

\title{Lagrangians for the Gopakumar--Vafa conjecture}

\author{Clifford Henry Taubes}
\givenname{Clifford Henry}
\surname{Taubes}
\address{Department of Mathematics\\
Harvard University\\\newline
Cambridge\\
Massachusetts 02138\\USA}
\email{chtaubes@math.harvard.edu}
\urladdr{}

\volumenumber{8}
\issuenumber{}
\publicationyear{2006}
\papernumber{3}
\startpage{73}
\endpage{95}

\doi{}
\MR{}
\Zbl{}


\received{22 January 2002}
\revised{}
\accepted{}
\published{22 April 2006 (in these proceedings)}
\publishedonline{22 April 2006 (in these proceedings)}
\proposed{}
\seconded{}
\corresponding{}
\editor{}
\version{}

\arxivreference{math.DG/0201219}

\dedicatory{Reproduced by kind permission of International Press
from:\newline 
{\rm Advances in Theoretical and Mathematical Physics, 
Volume 5 (2002) pages 139--163}}

\AtBeginDocument{\let\bar\wbar\let\tilde\wtilde\let\hat\what\def\notin{\not\in}}

\newtheorem{lemma}{Lemma}
\theoremstyle{definition}
\newtheorem{sstep}{Step}
\def\Step#1{\begin{sstep}\label{step#1}}
\def\NextStep#1{\end{sstep}\begin{sstep}\label{step#1}}
\newtheorem{ppart}{Part}
\def\Part#1{\begin{ppart}\label{part#1}}
\def\NextPart#1{\end{ppart}\begin{ppart}\label{part#1}}

\newcommand{\nc}{\newcommand}
\nc{\al}{\alpha}
\nc{\be}{\beta}
\nc{\de}{\delta}
\nc{\ga}{\gamma}
\nc{\La}{\Lambda}
\nc{\la}{\lambda}
\nc{\om}{\omega}
\nc{\si}{\sigma}
\nc{\Si}{\Sigma}
\nc{\ta}{\theta}
\nc{\va}{\varphi}
\nc{\ve}{\varepsilon}
\nc{\ze}{\zeta}
\nc{\ch}{{\mathcal H}}
\nc{\calO}{{\mathcal O}}
\nc{\pa}{\partial}
\nc{\we}{\wedge}
\nc{\bC}{\mathbb C}
\nc{\bP}{\mathbb P}
\nc{\bR}{\mathbb R}
\nc{\bZ}{\mathbb Z}
\nc{\iy}{\infty}
\nc{\op}{\oplus}
\nc{\ti}{\times}
\nc{\sub}{\subset}
\nc{\ra}{\rightarrow}

\begin{document}

\begin{asciiabstract}
This article explains how to construct immersed Lagrangian
submanifolds in C^2 that are asymptotic at large distance from the
origin to a given braid in the 3-sphere. The self-intersections of the
Lagrangians are related to the crossings of the braid. These
Lagrangians are then used to construct immersed Lagrangians in the
vector bundle O(-1) oplus O(-1) over the Riemann sphere which are
asymptotic at large distance from the zero section to braids.
\end{asciiabstract}

\begin{htmlabstract}
This article explains how to construct immersed Lagrangian submanifolds in
<b>C</b><sup>2</sup> that are asymptotic at large distance from the origin to a
given braid in the 3&ndash;sphere. The self-intersections of the Lagrangians
are related to the crossings of the braid. These Lagrangians are then
used to construct immersed Lagrangians in the vector bundle
O(-1)&oplus;O(-1) over the Riemann sphere which are asymptotic at large
distance from the zero section to braids.
\end{htmlabstract}

\begin{abstract}
This article explains how to construct immersed Lagrangian submanifolds in
$\mathbb{C}^2$ that are asymptotic at large distance from the origin to a
given braid in the 3--sphere. The self-intersections of the Lagrangians
are related to the crossings of the braid. These Lagrangians are then
used to construct immersed Lagrangians in the vector bundle $O(-1)
\oplus O(-1)$ over the Riemann sphere which are asymptotic at large
distance from the zero section to braids.
\end{abstract}

\maketitle


Gopakumar and Vafa \cite{GV} have have conjectured the existence of
a fundamental relationship between Gromov--Witten type invariants of
holomorphic curves in the vector bundle $O(-1) \op O(-1)$ over $\bP^1$ and
certain knot invariants, for example the Jones polynomial.  They came
to their conjecture by applying a fundamental observation of `t Hooft
\cite{H} in a string theoretic context on $T^*S^3$ 
described by Witten \cite{W}.  Subsequently, the scope of the
conjecture was expanded by Ooguri and Vafa \cite{OV}.  Successful
tests of the have been made, for example, by Labastida and Marino
\cite{LM}, Ramadevi and Sarkar \cite{RS}, Labastida, Marino and Vafa
\cite{LMV} and Aganagic, Klemm and Vafa \cite{AKV}.  
In the mean time, Faber and Pandarhapande \cite{FP}, Katz and Liu
\cite{KL} and Li and Song \cite{LS} have considered the mathematical
foundations for the conjecture and verified certain parts of it.

The verification of Gopakumar and Vafa's proposal has been slow, in
part because the string theoretic arguments have not provided a
geometric correspondence between a particular knot and a particular
set of holomorphic curves in $O(-1) \op O(-1)$.  Even so, it is a good
bet, verified in part by Katz and Liu \cite{KL}, Labastida, Marino and Vafa
\cite{LMV} and Aganagic, Klemm and Vafa \cite{AKV}, 
that such a correspondence
exists and that it is mediated by a suitable Lagrangian 3--manifold
sitting in $O(-1) \op O(-1)$.  To be specific, the knot should determine
the Lagrangian, and then a knot invariant should come as a suitable
count of compact, holomorphic curves with boundary on the Lagrangian.

This said, the mathematics of counting holomorphic curves with
boundary on a Lagrangian submanifold dates back to Floer's original
work on the Arnold conjecture \cite{F}.  Moreover, since Floer's work,
such counts have been considered by mathematicians in myriad
circumstances.  Yet, each new circumstance typically has new technical
problems to surmount; and in this regard, Katz and Liu \cite{KL} found the
story here to be typical.  

Curve counting theory aside, the proposed
mechanism via Lagrangians for Gopakumar and Vafa's conjecture
requires knots to provide Lagrangians in $O(-1) \op O(- 1)$.  This article
addresses the latter concern.  By way of preliminary remarks, note
that the construction described below produces a 2--dimensional
Lagrangian surface in $\bC^2$ from a knot in $S^3$, designed so that the
Lagrangian surface intersects all large radius 3--spheres as an
isotopy of the knot.  Moreover, if the starting knot is isotopic to
one that is mapped to itself via $S^3$'s antipodal map, then the
resulting Lagrangian is used to construct a 3--dimensional Lagrangian
in $O(-1) \op O(-1)$ that fibers over the equator in $\bP^1$ with the 
2--dimensional Lagrangian as fiber.  Thus, the construction below can be
viewed as one that constructs a Lagrangian in $O(-1) \op O(-1)$ from a
knot in $\bR\bP^3$.  

By the way, start with a knot in the unit radius sphere
in $S^3$ and Gromov's $h$--prinicple \cite{G} more or less asserts that there
is a Lagrangian surface in the unit ball of $\bC^2$ that intersects the
boundary 3--sphere as the given knot. This said, the construction
below provides explicit Lagrangians. In particular, a realization of
the knot as a braid provides a Lagrangian whose topology can be read
off directly from the properties of the braid.  

Here is how the remainder of this article is organized: \fullref{sec1}
 describes the construction of Lagrangians in $\bC^2$ from
braids in $S^3$ that intersect all large radius 3--spheres as a braid
that is braid isotopic to the original.  The initial steps construct
immersed disks, and it is then explained how the immersion points can
be smoothed to produce embedded, although not always orientable,
Lagrangian surfaces.  \fullref{sec2} explains how the double points of the
immersed disks relate to the crossings of a certain projection of the
original braid.  \fullref{sec3} relates the topology of three
Lagrangians coming from a triad of braids that arise in standard
discussions of skein relations.  \fullref{sec4} constructs Lagrangian
3--manifolds in $O(-1) \op O(-1)$ from certain Lagrangian surfaces in
$\bC^2$, and \fullref{sec5} views these $O(-1) \op O(-1)$
Lagrangians from a dual perspective on $T^*S^3$.  Note that the
discussion that follows owes much to conversations between the author
and Professor Cumrun Vafa to whom thanks is offered.

\section{The construction of Lagrangians}\label{sec1}

	My purpose here  is to describes a construction
that starts with a connected, $N$--stranded braid in $S^3$ and constructs of
an a properly immersed, Lagrangian disk in $\bC^2$ that intersects large
radiiall sufficiently large radius 3--spheres in $\bC^2$ as an $N$--stranded
braid that is braid isotopic to the original. The construction is
then generalized to obtain immersed Lagrangians for braids with more
than one component, and then generalized again to provide embedded
(but possibly non-orientable) Lagrangians.  The construction is
divided into ten steps. 3--spheres as a given $N$--stranded braid.  The
construction is generalized at the end to the case of links.

\Step 1
To start the construction, introduce $(z_1, z_2)$ to denote the standard, 
complex coordinates on $\bC^2,$ defined so that the symplectic form is given as
\begin{equation}
\om\equiv i 2^{-1} (dz1 \we d \bar{z}_1 + dz_2 \we d\bar{z}_2).\label{eq1}
\end{equation}
Next, introduce the `hyperk\"ahler' rotated complex coordinates 
\begin{equation}
a1 = 2^{-1/2} (z_1 - \bar{z}_2)\;\text{ and }\;   
a_2 = 2^{-1/2} (z_2 + \bar{z}_1) \label{eq2}
\end{equation}
with respect to which
\begin{equation}
\om = i 2^{-1}(da_1 \we da_2 - d \bar{a}_1 \we d \bar{a}_2).\label{eq3}
\end{equation}
Since this rotation is orthogonal, the metric in the $(a_1, a_2)$
coordinates is the standard one.  Note that the $a_1$--plane is a
Lagrangian plane in $\bC^2$ with respect to $\om$.  Of course, so is any
surface given as the zeros of a function of $(a_1, a_2)$ that is
holomorphic.  

In the subsequent steps, the complex coordinate $a_1$ is
written in terms of its real and imaginary parts as $a_1 = x^1 + i x^2$.
At the same time the coordinate $a_2$ is written somewhat perversely as
$a_2 = p_2 + i p_1$.  In terms of these real coordinates,
\begin{equation}
\om = dp_1 \we dx_1 + dp_2 \we dx^2 .\label{eq4}
\end{equation}
The choice of coordinates here is meant to stress an implicit
identification below between the $a_2$--direction in $\bC^2$ and the fiber of
the cotangent bundle of the $a_1$--plane.  More to the point, this gives a
symplectic identification between $\bC^2$ and the cotangent bundle, 
$T^*\bC$, of the $a_1$--plane.
	
The reason for making such an identification is as follows: If $f$ is
any smooth, locally defined function of the coordinates $(x_1, x_2)$, then
the graph of $df$ defines a Lagrangian surface 
in $T^*\bC$ and hence in $\bC^2$.
To be explicit, the locus of points in $T^*\bC$ where 
$(x^1, x^2, p_1 = \pa_1f, p_2 = \pa_2f)$ is a Lagrangian surface.

\NextStep 2 Fix attention on some given connected, $N$--stranded
braid $K$.  To be more precise about what this means here, identify
$S^1$ with the unit circle in $\bC$ with coordinate $\ze$ such that
$|\ze| = 1$.  This done, then $K$ can be viewed as an embedded circle
in $S^1\ti \bC \sub \bC^2$ obtained as the image of a map from $S^1$.
In particular, such a map should send $\eta\in S^1$ to $(\zeta=\eta^N,
z=\ga(\eta))$ with $\ga$ a map to $\bC$ that separates pairs of points
in $S^1$ whose ratio is a non-trivial, $N$'th root of unity.  In this
regard, it proves convenient below to write $\eta= e^{i\ta/N}$ with
$\ta\in [0, 2\pi N]$.  Thus, $\ga$ is a function of $\ta$ that is
periodic with period $2\pi N$ and $\ze = e^{i\ta}$ is a function of
$\ta$ with period $2\pi$.  This notation identifies the unit circle in
the $a_1$--plane to the interval $[0, 2\pi]$ with its endpoints
identified.

Note that the graph of any $2\pi N$ periodic complex function, $\ga$,
of $\ta$ defines a braid as long as the $N$ values $\{\ga(\ta + 2\pi
k)\}_{0\leq k<N}$ are distinct at each point $\ta\in [0, 2\pi]$.  In
particular, this last point of view will be taken here.
	
For example, a trivial $N$--stranded braid sits in $S^1\ti \bC$ as the
graph of $2\pi N$ periodic map $\ga\co  S^1 \ra \bC$ given by
\begin{equation}
\ga(\ta) = R^{-1} e^{i\ta/N} , \label{eq5}
\end{equation}
where $R > 0$ is any constant.

	Reference has been made at the outset to braids in $S^3$.  Take
this to mean the following: View $S^3$ as the unit radius sphere about
the origin in $\bC^2$ and identify $S^1\ti \bC$ 
with its image in $S^3$ via the embedding that sends a pair 
$(\ze, z)$ to $(\ze, z)/(1 + |z|^2)^{1/2} \sub \bC^2$.  This
done, an $N$--stranded braid in $S^3$ signifies the image of $S^3$ of such a
braid in $S^1\ti \bC$.
	
Two braids in $S^3$ are said below to be `braid isotopic' if they
are isotopic through a 1--parameter family of braids.

\NextStep 3 This step constructs an embedded, Lagrangian cylinder
in a neighborhood of $S^1\ti \bC$ in $\bC \ti \bC$ that intersects
$S^1\ti \bC$ as the given braid $K$.  For this purpose, introduce the
complex function $\ga$ on $S^1$ that defines $K$ and write $a_2 =
\ga_1 + i \ga_2$.  With the comments at the end of the preceding step
in mind, $K$ can be written as the intersection of $S^1 \ti \bC \sub
\bC^2$ with a Lagrangian cylinder defined in a neighborhood in $\bC^2$
of $S^1\ti \bC$ provided that the following is true: The section
$\ga_2 dx^1 + \ga_1 dx^2$ of $T^*\bC|_{S^1} = S^1 \ti \bC$ extends to
a section of $T^*\bC$ over a cylindrical neighborhood of the unit
circle in the $a_1$--plane as the differential of a function that is
$2\pi N$ periodic on the constant radius circles.  Thus, the goal is
to find a function on the $a_1$--plane, $2\pi N$ periodic on constant $r$
circles, whose partial derivative in the $x^1$--direction restricts to the
unit circle as $\ga_2$ and whose partial derivative in the $x^2$ direction
restricts to the unit circle as $\ga_1$.
	
To find such a function, it proves useful to introduce the radial
coordinates $r \geq 0$ and $\ta$ for the $a_1$--plane and write $x^1 =
r \cos\ta$ and $x^2 = r \sin\ta $.  This done, then
\begin{equation}
\ga_2 dx^1 + \ga_1 dx^2 = (\ga_2 \cos\ta + \ga_1 \sin\ta) dr 
+ (\ga_1 \cos\ta - \ga_2 \sin\ta) d\ta ,\label{eq6}
\end{equation}
and the task at hand is to find a function, $f$, of $r$ and $\ta$ such that
\begin{align}
\bullet &\quad f(r, \ta + 2\pi N) = f(r, \ta).\label{eq7}\\
\bullet &\quad \pa_rf |_{r=1} = \ga_2 \cos \ta + \ga_1 \sin \ta.\nonumber\\
\bullet &\quad \ta_{\ta}f|_{r=1} = \ga_1 \cos \ta - \ga_2 \sin \ta.\nonumber
\end{align}
	
There is one immediate requirement for $f$'s existence, which is that
\begin{equation}
\int^{2\pi N}_0(\ga_1(\ta) \cos\ta 
- \ga_2 \sin \ta) d\ta = 0\label{eq8}
\end{equation}
since this integral is meant to be $f(1, 0) - f(1, 2\pi N)$.  In this
regard, notice that any given $2\pi N$ 
periodic map $\ga$ to $\bC$ with $\{\ga(\ta +2\pi k)\}_{0\leq k<N}$ 
distinct at all values of $\ta$ can be homotoped through
such maps to one that obeys \eqref{eq8}.  In particular, such a homotopy does
not change the braid isotopy class of the corresponding braid.
Indeed, if $\ga\co  S^1 \ra \bC$ represents an 
$N$ stranded braid and if $c \in\bR$, then
$\ga'\equiv\ga + c e^{-i\ta}$ also has $N$ 
distinct values at each point for
$\{\ga'(\ta + 2\pi k)\}_{0\leq k<N}$.  Meanwhile, the value of the $\ga'$
version of \eqref{eq8} differs from the value of the $\ga$ version by $2\pi N c$ so
there is a unique such c for which the $\ga'$ version of \eqref{eq8} is zero.
This understood, agree henceforth to restrict attention to those 
maps $\ga$ where \eqref{eq8} holds.  
	
Given that \eqref{eq8} holds, then there exists a bonafide, $2\pi N$ periodic
function $\si$, of $\ta$ whose partial with respect 
to $\ta$ is equal to $\ga_1 \cos \ta- \ga_2 sin \ta$.  
This understood, then
\begin{equation}
 h \equiv (r - 1) (\ga_2 \cos \ta + \ga_1 \sin \ta) + \si(\ta) .\label{eq9}
\end{equation}
satisfies the conditions in \eqref{eq7} and so the graph of $dh$ in $T^*\bC = \bC^2$
provides an example of the required Lagrangian, at least near the unit
circle in the $a_1$--plane.  As demonstrated in the next step, the
Lagrangian defined by \eqref{eq9} is per force embedded near this circle, but
perhaps not everywhere.

\NextStep 4 As remarked at the end of the previous step, the
Lagrangian cylinder defined by the graph of the differential of the
function in \eqref{eq9} may have immersion points where $|a_1|$ differs
substantially from 1.  This step and the next describe how to define a
properly embedded, Lagrangian cylinder, defined near the $|a_1| = 1$
circle and where $|a_1| \geq 1$ that intersects every constant $|a_1|$ slice
as a braid that is isotopic to the original.
	
To start this construction, represent the given braid using, as
described, a $2\pi N$ periodic map 
$\ga = \ga_1 + i \ga_2\co  S^1 \ra \bC$ with distinct
values for $\{\ga(\ta + 2\pi k)\}_{0\leq k<N}$ at all points.  By way of
shorthand, introduce $\al\equiv \ga_2 \cos \ta + \ga_1 \sin \ta$
and $\be\equiv 1 \cos \ta - \ga_2 \sin\ta$.  
Note that the pair $(\al, \be)$ are $2\pi N$ periodic, and the $N$ pairs
$\{(\al(\ta+ 2\pi k), \be(\ta + 2\pi k))\}_{0\leq k<N}$ 
are distinct at each $\ta$ if
and only if such is the case for $\{\ga(\ta + 2\pi k)\}_{0\leq k<N}.$
	
Reintroduce the function, $h$, in \eqref{eq9}; its partial derivatives determine the 
Lagrangian cylinder from the preceding step.  In particular, these
derivatives are 
\begin{align}
\bullet &\quad \pa_rh = \al .\label{eq10}\\
\bullet &\quad \pa_{\ta}h = (r - 1) \pa_{\ta} + \be .\nonumber
\end{align}
By virtue of continuity and the fact that $\{(\al(\ta+ 2\pi k),
\be(\ta + 2\pi k))\}_{0\leq k<N}$ has $N$ distinct pairs at each
$\ta$, the differential $dh = \pa_rh dr + \pa_{\ta}h d\ta$ has the property 
that $\{dh|_{\ta+2\pi k}\}_{0\leq k<B}$ also has $N$ distinct values at
each point of the constant $r$ circle if $|r - 1|$ is not too big.  In
particular, there exists some $\de > 0$ for which such is the case
when $|r - 1| \leq 2 \de$; and this implies that the graph of $dh$
defines an embedded Lagrangian cylinder where $|r - 1| \leq 2 \de$.
Of course an upper bound for $\de$ is determined by the braid map
$\ga$, but there is no positive lower bound to the choice of $\de$ to
use here and in the subsequent discussions.  In particular, the
condition $\de < 10^{-3}$ is implicitly enforced.
	
Here is a reformulation of this last point for use below: As long as
$|s| \leq 2 \de$, then the graph over $S^1$ 
in $T^*\bC$ of the 1--form $\al dr + (s \pa_{\ta}\al + \be) d\ta$
is such that its values at the points $\{\ta + 2\pi k\}_{0\leq k<N}$ 
are distinct at each $\ta$ and so defines a braid that is braid
isotopic to the original (that with $s = 0$).  More generally, as long
as $\ve$ is not zero and $|s| <2\de $ then the 1--form
\begin{equation}
\ve\al dr + (s \pa_{\ta}\al + \be) d\ta \label{eq11}
\end{equation}
also has this same property.  Thus, the graph over the circle of the
1--form in \eqref{eq11} defines a braid that is braid isotopic to the original
braid defined by $\ga$.

\NextStep 5
This step uses the observation in \eqref{eq11} to obtain the promised cylinder 
from \fullref{step4}.  For this purpose, replace the function in \eqref{eq9} and \eqref{eq10} by
\begin{equation}
f = \de (r - 1) (r + \de)^{-1}\al + \si .\label{eq12}
\end{equation}
where $\si$ is as before, $\pa_{\ta}\si = \be$.  This choice gives
\begin{align}
\bullet &\quad  \pa_rf = \de (1 + \de) (r + \de)^{-2} \al .\label{eq13}\\
\bullet &\quad  \pa_{\ta}f=\de(r-1)(r+\de)^{-1}\pa_{\ta}\al+\be.\nonumber
\end{align}
As before, the graph of df defines a Lagrangian.  In particular, with
\eqref{eq11} in mind, it follows that $df$ on any $r \geq 1 - 2\de$ circle has
distinct values at $\{\ta + 2\pi k\}_{0\leq k<N}$ 
for each $\ta$ and so the graph
of $df$ is a properly embedded, Lagrangian cylinder in the $|a_1| \geq 1
- 2 \de$ portion of $\bC^2$ whose intersection with any constant $|a_1|$ slice
is a braid that is braid isotopic to the original.

\NextStep 6 This step explains how to extend the cylinder defined
in the previous step to the $|a_1| \leq 1 - 2 \de$ portion of $\bC^2$ capping
the $r = 1 -\de$ slice of this cylinder with a closed, immersed,
Lagrangian disk in the $r \leq 1 - \de$ of $\bC^2$.  The self intersection
points of this added disk are described in a subsequent step.

The construction of this extension starts by returning to the example
of the trivial $N$ stranded braid where 
$\ga(\ta) = e^{i\ta/N}$.  This braid has 
$\ga_1= \cos(\ta/N)$ and $\ga_2 = \sin(\ta/N)$ 
and so \eqref{eq8} is satisfied.  Of course,
there is a Lagrangian that extends this particular braid, it given by
the locus of points $(a_1, a_1^{1/N})$ 
in $\bC^2$ with $a_1 = r e^{i\ta}$.  This extension
is given as the graph of $df_N$, where
\begin{equation}
f_N =  (1 + 1/N)^{-1} r^{1+1/N} \sin((1 + 1/N) \ta).\label{eq14}
\end{equation}
With $f_N$ understood, the differential of any function that
interpolates between $f$ in \eqref{eq12} where $r \geq 1 - \de$ and $f_N$ in \eqref{eq14} where
$r$ is near zero defines an immersed Lagrangian disk with the requisite
properties.  For example,
\begin{equation}
(1 - \chi) f + \chi f_N \label{eq15}
\end{equation}
is such an interpolating function with $\chi$ any function of $r$ that
equals 1 near $r = 0$ and 1  where $r \geq 3/4$.

\NextStep 7
The double points of the Lagrangian defined by \eqref{eq15} can be related 
directly to properties of the original braid $\ga$.  
These relations are described below in 
\fullref{sec2} for a more sophisticated version of the function $f$ that
appears in \eqref{eq15}.  The  
description of this new $f$ requires the specification of a small and
positive constant $\de$.   
Given $\de$, fix a smooth function, $\chi_{\de}$, of the coordinate r
that has the following properties: 
\begin{align}
\bullet&\quad 
\hbox{$\chi_{\de} = 1$  where   $r \leq 1 - 2\de$.}\label{eq16}\\ 
\bullet&\quad 
\hbox{$\chi_{\de} = 0$  where   $r \geq 1 - \de$.}\nonumber\\ 
\bullet&\quad 
\hbox{$\chi_{\de} = \de^{-1} (1 - \de - r)$  where   
$1 - 2 \de + \de^2 \leq r \leq 1 -\de - \de^2$.}\nonumber\\
\bullet&\quad 
\hbox{$\pa_r\chi_{\de} \leq 0$.}\nonumber
\end{align}
Further we require:
\begin{itemize}
\item[$\bullet$] 
Set  $r^* \equiv 1 - 2 \de + 2 \de^4$  and require that $r^*$ be the unique 
value of $r$ where $\chi_{\de}r^{1+1/N}$  achieves its maximum, 
and require that this maximum be non-degenerate in the sense 
that 
\begin{itemize}
\item[(a)]  	
$\pa_r(\chi_{\de} r^{1+1/N}) > 0$ where $r < r^*$.
\item[b)]
$\pa_r(\chi_{\de} r^{1+1/N}) < 0$ where $r > r^*$.
\item[c)]
$\pa_r(\chi_{\de} r^{1+1/N}) = r^*-r$  where  
$r^* - \de^4 < r < r^* + \de^4$.
\end{itemize}
\item[$\bullet$]
Where  $1 - \de+ \de^2 < r \leq 1 - \de$,
\begin{itemize}
\item[a)]
$\chi_{\de}|\pa_r\chi_{\de}|^{-1} \leq 10 \de^2$. 
\item[b)]
$|\pa_r(\chi_{\de}r^{1+1/N})|$   is decreasing.
\end{itemize}
\item[$\bullet$] 
$\pa_r^2\chi_{\de}>0$ and $|\pa_r\chi_{\de}|\leq 100\de|\pa_r^2\chi_{\de}|$  
where  $1 - \de + \de^2/2 \leq r \leq 1 - \de$.  
\end{itemize}

Note that the third to last point above asks only that $\chi_{\de}
r^{1+1/N}$ behave in a uniformly quadratic fashion near its maximizer,
$r^*$.  Meanwhile, the final point two points can be achieved by
requiring $\chi_{\de}$ to vanish as $r \ra 1 -\de $ as a multiple of
the exponential of the function $-(1 - \de - r)^{-2}$.

Fix a second smooth function, $\chi$, of r that has value 1 where $r
\leq 1/2$, value 0 where $r \geq 3/4$ and whose derivative is nowhere
greater than 8.

With the preceding understood, replace the function $f$ in \eqref{eq15}, by
\begin{equation}
f\equiv (1 - \chi) f + \chi_{\de} f_N .\label{eq17}
\end{equation}
By construction, the graph of $df$ then defines a smooth, properly
immersed Lagrangian disk, $L$, in $\bC^2$ whose $|a_1| > 1 -\de $
portion is embedded and intersects every $|a_1| \geq 1 - \de$ slice of
$\bC^2$ transversely as a braid that is braid isotopic to the original
one.

In addition, if $\ga$ is replaced by $\ve\ga$ with $\ve > 0$ and very
small, (so representing an isotopy of the original braid to one with
distance $\calO(\ve)$ from the $a_1$--plane), then the graph of the
differential of the $\ve\ga$ version of $f$ in \eqref{eq17} produces a
properly immersed, Lagrangian disk in $\bC^2$ that is embedded near
the $|a_1| = 1$ slice, embedded where $|a_1|\geq 1$, 
and intersects every radius 1 or larger 3--sphere transversely
and in a braid that is braid isotopic to the original.
	
\NextStep 8
Although the Lagrangian defined by the differential of
the function in \eqref{eq17} has various virtues, it may not be the most
useful for certain applications.  This step describes a second
Lagrangian in $\bC^2$ with a somewhat different suite of properties.  In
particular, the construction here facilitates comparisons when
non-isotopic braids differ by a strand crossing.  However, the down
side here is that the Lagrangians from this step may only intersect
all sufficiently large radius spheres as a braid isotopy of the
original braid.

To start the construction, choose, as before, a function $\chi$ of the
radial coordinate $r$, where now $\chi$ can have value 1 near $r = 0$ and value
0 at large $r$.  Let $\ga$ be a given braid 
and again introduce $\al, \be$ and $\si$.
With $f_N$ as in \eqref{eq14} fix some $\de > 0$ to define
\begin{equation}
f_{\bullet} (1 - \chi) (- r^{-1} \al + \si) +\chi\de f_N .\label{eq18}
\end{equation}
Note that where $r$ is large and so  $\chi= 0$, 
the differential of $f_{\bullet}$ is given by
\begin{equation}
df_* = r^{-2} \al dr + (- r^{-1} \pa_{\ta}\al + \be) d\ta .\label{eq19}
\end{equation}
The advertised new Lagrangian is defined by the graph of $df_*$. Note
that the discussion in \fullref{step4}'s final paragraph justifies the claim
that this new Lagrangian intersects all sufficiently large radius
spheres transversely as a braid isotopy of the original braid.

\NextStep 9 This step constructs Lagrangians in $\bC^2$ that
intersect the large radius 3--spheres as a braid isotopy of a given $N$
stranded, but multiple component braid.  In particular, after suitably
parametrizing the braid, the construction is essentially identical to
that described in the previous steps.  To start, suppose that the
braid has some $n$ components, $\{\ga_1, \ga_2, \ldots , \ga_n\}$
where each $\ga_j$ is a function of $\ta$ that is periodic with period
$2\pi N_j$.  Here, $\Si_j N_j = N$.  Now, take the parameter $\de$ to
be very small, and for each $j$, use the chosen $\de$ to construct
that $\ga_j$ version of the function in either \eqref{eq15}, \eqref{eq17}. Call it
$f_j$ and let $L_j \sub \bC^2$ denote the corresponding Lagrangian.
Then, the claim is that $L = \cap_j L_j$ is an immersed Lagrangian with
the desired properties.

The proof of this claim requires only a verification that $L$'s
intersection with all large radius 3--spheres is braid isotopic to the
original braid.  For this purpose, note that no $\ta\in [0, 2\pi]$
exist where a pair from any $\ga_i$ version of $\{\al(\ta + 2\pi k),
\be(\ta + 2\pi k)\}_{0\leq k<N}$ coincides with one from the $\ga_j$
version when $i \neq j$.  This understood, it follows by continuity
from \eqref{eq11} that choosing $\de > 0$ and small guarantees that the
corresponding $df_i$ and $df_j$ have disjoint graphs where $r \geq 1 -
\de$.  Thus, $L_i$ and $L_j$ are disjoint where $r \geq 1$ and their
intersection with any $r \geq 1 - \de$ slice of $\bC^2$ is a braid
that is braid isotopic to the original.

An alternate construction takes each $f_j$ to be the $\ga_j$ version
of the function defined by \eqref{eq19}, and then takes $L_j$ to be the
corresponding immersed, Lagrangian disk.  This understood, set $L
\equiv\cup_j L_j$.  The latter is immersed, and as follows from \eqref{eq19}
using perturbation theory, it intersects all sufficiently large radius
3--spheres transversely in a braid that is braid isotopic to the
original.

\NextStep {10} This step describes how to modify an immersed
Lagrangian on some small neighborhood of its immersion points to
obtain an embedded, but higher genus Lagrangian.  In this regard, note
that the resulting Lagrangian may not be orientable.  In particular,
this situation occurs when the initial Lagrangian has
self-intersection points with positive local degree.  By the way, the
existence of such a modification has surely been known for years by
experts, but as the construction is relatively straightforward, it is
worth relating the details.

The first point to make is that any immersed Lagrangian can be
modified on any given neighborhood of its singular points (with out
changing the genus) so that the result has only transversal, double
point self-intersections.  Moreover, if the original singular set is
compact, then this modification produces only a finite set of such
intersections.  The definition of such a modification exploits the
4--dimensional version of the following basic and well known lemma:

\begin{lemma}\label{lem1}  
Let $n \geq 1$ be an integer, $X$ be a $2n$--dimensional manifold with
a symplectic form, and $L \sub X$ be an $n$--dimensional Lagrangian
submanifold.  Then each point of $L$ has a neighborhood with
coordinates $(x^1,\ldots,x^n, p_1,\ldots , p_n)$ in which the $\{p_j =
0\}_{1\leq j\leq n}$ slice is $L$ and to which the symplectic form
restricts as $\Si_j dp_j dx^j$.
\end{lemma}

Given the lemma, a first perturbation of $L$ puts a neighborhood of
some given singular point in the appropriate form.  For this purpose,
select a singular point $z \in L$ and take a coordinate system
centered on this point as described by the lemma.  A perturbation of $L$
near $z$ is defined by the locus where $\{p_j = \pa_jh\}$ where $h$ is a
function on $L$ whose partial derivatives are small.  In particular, as
there are functions defined near the origin in $\bR^n$ with any given
vector as differential at 0 and any given symmetric matrix as Hessian,
Sard's theorem provides perturbations of $L$ that stay arbitrarily
close to $L$, contains $z$ and result in a new Lagrangian, $L'$, with a
transversal and purely double point self-intersection at $z$.

This local construction understood, a straightforward extension
produces a Lagrangian with all singular points as desired.  The
details of the extensions are tedious and omitted.

Given that a singular point in the Lagrangian $L$ is isolated and a
transversal double point, the modification to make a Lagrangian that
is embedded with one less singular point procedes as follows: First,
fix one of the sheets of the Lagrangian on a neighborhood of a
singular point and introduce the Lemma's coordinates with center on
the singular point.  This done, the other sheet can be perturbed
without introducing new singular points so that it intersects the
original in the origin of these coordinates and so that a neighborhood
of the origin in this sheet coincides with the locus where $x^1 = x^2 =0$.
This understood, introduce the complex coordinates $a_1 = x^1 + i x^2$ and
$a_2 = p_2 + i p_1$ with respect to which the symplectic form is given by
\eqref{eq3} and $L$'s intersection with a neighborhood of the origin is the
locus where $a_1a_2 = 0$.

Now, consider the perturbation of $L$ in this neighborhood given by the
locus where $a_1a_2 = i\ve$ with $\ve$ some non-zero, small positive constant.
The latter locus is a smooth, Lagrangian submanifold.  Moreover, if $\ve$
is small, then its intersection with the complement of a small radius
ball about the origin consists of two annuli, one very close to the
$a_1$--plane and the other close to the $a_2$--plane.  Here, the annulus that
near the $a_1$--plane is the locus of points where $p_1 = \ve\pa_1 \ln(r)$ and
$p_2 = \ve\pa_2 \ln(r)$ where $r = (x_1^2 + x_2^2)^{1/2}$.  
Meanwhile, the annulus
that is close to the $a_2$--plane is defined by the analogous locus where
the role of the pair $(x^1, x^2)$ is switched with that of $(p_1, p_2)$.

These last points understood, fix a non-increasing function $\chi$ on $[0,
\iy)$ that equals 1 on $[0, 1]$, vanishes on $[3, \iy)$ and whose
derivative is no larger in absolute value than 1.  This done, replace
the annulus close to the $a_1$--plane by the locus of points where
$p_1$ and $p_2$ are the respective partial derivatives of
$\ve\chi(r/\ve) \ln(r)$.  At the same time, replace the annulus close
to the $a_2$--plane in the analogous fashion.  The result is, for small
$\ve$, a Lagrangian that has one less self-intersection point than the
original and agrees with the original in the complement of a small
neighborhood of the chosen self-intersection point.

By the way, this construction respects the given orientation on the
two intersecting sheets of the original Lagrangian only when the local
intersection number of the two sheets is $+1$.  For topological reasons,
it is impossible to remove a local intersection with intersection
number $-1$ using a local modification that preserves the orientations
on the intersecting sheets.\end{sstep}

\section{Immersion double points and crossings}\label{sec2}

The self intersection points of the Lagrangians defined from
either \eqref{eq15}, \eqref{eq17} or \eqref{eq19} can be directly related to properties of the
original braid.  This is done here for the small $\de$  versions of the
Lagrangian given by \eqref{eq17}.

To start, suppose that the constant $\de$ that appears in \eqref{eq12} and
\eqref{eq17} is taken very small (remember that the discussion in \fullref{sec1} is
valid as long as $\de$ is positive no matter how small).  The fact is
that the integer $N$ determines an upper bound for the application
that follows, but such an upper bound is not explicitly derived.  With
$\de$ small, replace the map $\ga$ that defines a given braid by
$\ve\ga$ with $\ve > 0$ and very small.  This done, then the
self-intersection points of the Lagrangian from the resulting $f$ in
\eqref{eq17} can be interpreted in terms of the crossings of the original
braid.  The purpose of this section is to explain how this comes
about.  The discussion that follows is divided into six parts.

\Part 1  
To begin the story, remark that with $\ga$ given and the function $f$ defined 
from $\ga$ as in \eqref{eq12} then
\begin{equation}
f_{\ve}\equiv (1 - \chi) \ve f +\chi_{\de} fN .\label{eq20}
\end{equation}
is the $\ve\ga$ version of \eqref{eq17}.  Now, note that $f_{\ve} = f_N$ where
$r \leq 1/2$ and so the Lagrangian that is defined by $df_{\ve}$
intersects the $|a_1| < 1/2$ portion of $\bC^2$ as an embedded disk.
When $\ve$ is small, such is the case where $|a_1| < 1 - 2 \de$ for
the following reason: Where $r \leq (1 - 2 \de)$, the differential of
$f_{\ve}$ differs from that of $\ve$, $f_N$ by a term no larger than
$\ve|df|$, and so for small the $|a_1| \leq (1 - 2 \de)$ portion of
the Lagrangian defined by $df_{\ve}$ is a small perturbation of that
defined by $df_N$.  In particular, as the latter is embedded, so the
$|a_1| < 1 - 2 \de$ portion of the former is also.

With the preceding understood, it follows that the immersion points of the 
Lagrangian in question all lie where $1 - 2 \de \leq |a_1| \leq 1$.  To
study these points, note first  
that the function $f_{\ve}$ where $1 - 2 \de\leq r \leq 1$ is given by
$f_{\ve} = \ve f + \chi_{\de} f_N$, and so self-intersection 
points in the graph of $df_{\ve}$ are at $(r, \ta)$ where the
values at $\{(r, \ta + 2\pi k)\}_{0\leq k<N}$  
of the $2\pi N$ periodic 1--form 
\begin{equation}
[\ve\pa_rf + \mu^{-1}\pa_r(\chi_{\de} r^{\mu}) \sin(\mu\ta)] dr 
+ [\ve\pa_{\ta}f + \chi_{\de} r^{\mu} \cos(\mu\ta)] d\ta \label{eq21}
\end{equation}
are not pairwise distinct.  To investigate where these $(r, \ta)$
occur, it proves useful to separate the search into three regimes.
The first occurs where $1 - 2 \de + \de^2 \leq r \leq 1 - \de -
\de^2$, the second where $1 - 2\de < r < 1 - 2 \de + \de^2$ and the
third where $1 - \de - \de^2 < r < 1 - \de$.
 
\NextPart2 
In the first regime, \eqref{eq16} implies that the 1--form in \eqref{eq21} looks like
\begin{equation}
\de^{-1} [\mu^{-1} \sin(\mu\ta) + \calO(\ve + \de)] dr 
+ [\de^{-1} (1 - 2\de + r) \cos(\mu\ta) + \calO(\ve)] d\ta .
\label{eq22}
\end{equation}
In particular, as the values at $\{\ta + 2\pi k\}_{0\leq k<N}$ of the
pair $(\cos(\mu\ta), \sin(\mu\ta))$ define a set of $N$ distinct
elements for each $\ta$, so do the values of the form in \eqref{eq19} at
$\{(r, \ta + 2\pi k)\}_{0\leq k<N}$ when both $\de$ and $\ve$ are
small.

\NextPart3   
Consider next the second regime, that where $1 - 2 \de\leq r \leq 1
- 2 \de + \de^2$.  Here, $\chi_{\de}$  
is close to 1 but the deriviative of $\chi_{\de}r^{\mu}$ has
a zero so the 1--form in \eqref{eq21} appears schematically as
\begin{equation}
[\ve\de_rf + \mu^{-1}\pa_r(\chi_{\de} r^{\mu} \sin(\mu\ta)] dr 
+  [cos(\mu\ta) + \calO(\ve)] d\ta.  \label{eq23}
\end{equation}
Now, given \eqref{eq23} and small $\ve$, perturbation theory precludes less
than $N$ distinct elements in the set of values of \eqref{eq23} at a given
$\{(r, \ta + 2\pi k)\}_{0\leq k<N}$ unless $\ta$ is close to a point
where $\{\cos(\mu\ta + 2\pi k/N)\}_{0\leq k<N}$ has less than $N$
distinct elements.  In this regard, a glance at the graph of the
cosine function indicates that there are $2(N - 1)$ points in $[0,
2\pi]$ where $\{\cos(\mu\ta + 2\pi k)\}_{0\leq kNN}$ has less than $N$
distinct elements and at such points, this set has precisely $N - 1$
distinct elements.  Moreover, the coincidence of a pair of elements of
this set at these special $\ta$ points is achieved in a manner that is
non-degenerate in the following sense: If $\ta_*\in [0, 2\pi]$ is one
of these special points, and if $k_*\in \{1,\ldots, 2\pi N\}$ is such
that $\cos(\mu\ta_* + 2\pi k_*/N) = \cos(\mu\ta_*)$, then the
derivative of $\cos(\mu\ta + 2\pi k_*/N) - \cos(\mu\ta)$ at $\ta_*$ is
non-zero.

By the way, this count of points in $[0, 2\pi]$ where $\{\cos(\mu\ta +
2\pi k/N)\}_{0\leq k<N}$ has less than $N$ distinct elements arises
from the fact that each $k\in \{1,\ldots , N - 1\}$ determines
precisely two values for $\ta\in [0, 2\pi]$ where $\cos(\mu\ta)$ and
$\cos(\mu\ta + 2\pi k/N)$ are equal.

These last remarks understood, an application of perturbation theory
finds, for small $\ve$, precisely $2(N- 1$) points $\ta\in [0, 2\pi]$
for which the set of values of the $d\ta$ component in \eqref{eq20} at $\{(\ta
+ 2\pi k)\}_{0\leq k<N}$ has less than $N$ distinct elements, and at
such a point, this set then has precisely $N - 1$ elements.  Moreover,
each such point in $[0, 2\pi]$ will be very close (for small $\ve$) to
a point where there are fewer than $N$ distinct elements in the set of
values of $\cos(\mu\ta)$ at $\{\ta + 2\pi k\}_{0\leq k<N}$.  Use
$\La\sub [0, 2\pi]$ denote the set of those $\ta$ where the $d\ta$
component of \eqref{eq20} at $\{\ta + 2\pi k\}_{0\leq k<N}$ has fewer than $N$
elements.

As demonstrated by a second application of perturbation theory, the
fifth point of \eqref{eq16} has the following implication: Given that $\de$ is
small and then $\ve$ is very small, each $\ta_*\in\La$ is the
$\ta$--component of a unique point $(r_*,\ta_*)$ with $r_*\in [r^* -
\de^4, r^* + \de^4]$ where the set of values of the whole 1--form in
\eqref{eq23} at $\{(r_*, \ta_* + 2\pi k)\}_{0\leq k<N}$ has less than $N$ (and
thus $N - 1$) elements.

Given all of the above, then it follows that the Lagragian in $\bC^2$
defined by the differential of the small $\de$ and very small $\ve$
version of $f_{\ve}$ in \eqref{eq20} has precisely $2(N - 1)$ double points
where $|a_1|$ lies between $1 -\de + \de^2$ and $1 -2\de$.  Moreover,
the arguments ust given establish that each of these self intersection
points of the Lagrangian is transversal.  Meanwhile, the discussion
below in \eqref{part6} explains why these self intersection points all
contribute the same local sign to any count of a self intersection
number of $L$.

\NextPart4
Consider now the third regime, that where $1 - \de -
\de^2 \leq r < 1 - \de$.  Here, it proves useful to break this regime
into two parts, the first where $\chi_{\de}\geq\de\ve$ and the second
where this last condition does not hold.  In this first regime,
$\pa_r\chi_{\de} = - 10^{-1}\de^{-1} \ve w$ where the function w is
greater than 1 by virtue of the second to last point in \eqref{eq16}.  Thus,
\eqref{eq21} has the schematic form
\begin{equation}
10^{-1}\de^{-1}\ve [- w\mu^{-1}\sin(\mu\ta) + \calO(\de)] dr 
+ [\ve\pa_{\ta}f + \chi_{\de} r^{\mu} \cos(\mu\ta)] d\ta .
\label{eq24}
\end{equation}
In particular, when $\de$ is small, then the values of the $dr$
component of \eqref{eq24} at the points in $\{(r, \ta + 2\pi k)\}_{0\leq k<N}$
is a set of fewer than $N$ distinct elements provided that $\ta$ is
near one of the $2(N - 1)$ points where $\{(\sin(\mu\ta+ 2\pi
k)\}_{0\leq kNN}$ has fewer than $N$ distinct elements.

Now, not all of these $2(N - 1)$ points in $[0, 2\pi]$ correspond to
self-intersection points of the Lagrangian with $r$ in the prescribed
range.  Indeed, when $\ve$ and $\de$ are very small, 
then the form of the $d\ta$ component of \eqref{eq24} 
forces a self-intersection point at $(r, \ta)$ with $\ta$
very close to such a $\ta_*$ and with $r$ in the prescribed range provided
that the following requirment is met: The expressions
\begin{equation}
\be(\ta_* + 2\pi k) - \be(\ta_* + 2\pi k')\quad   
\text{and}\quad \cos(\mu\ta_* + 2\pi k/N) - \cos(\mu\ta_* + 2\pi k'/N)
\label{eq25}
\end{equation}
have opposite sign when $k \neq k'\{0,\ldots, N - 1\}$ 
are chosen to make $\sin(\mu\ta_* + 2\pi k) = 
\sin(\mu\ta_* + 2\pi k')$.  
Moreover, if the left most difference in \eqref{eq25} is non-zero for all of
the $2(N - 1)$ possibilities for $\ta_*$, then perturbation theory
guarantees a 1--1 correspondence between the self intersection points
in the third regime and those $\ta_*$ where the just stated
requirement is met.  This guarantee also comes with a rider to insure
that these self intersection points are all transverse double points.

\NextPart5   
Consider the final part of the third regime where 
$\chi_{\de}<\de\ve $.  Here, \eqref{eq21} has the schematic form
\begin{equation}
[\ve\de(1 + \de) (r + \de)^{-2} \al
+ \mu^{-1}\pa_r(\chi_{\de} r^{\mu} \sin(\mu\ta)] dr 
+ \ve[\be + \calO(\de)] d\ta .
\label{eq26}
\end{equation}
In this regard, note that when $\de$ is small, then the $d\ta$
component of \eqref{eq26} has $N$ distinct values except possibly near points
in $[0, 2\pi]$ where the function $\be$ has fewer than $N$ distinct
values.

Introduce the term `twisted crossing point' to denote a point
$\ta_*\in [0, 2\pi]$ where the set $\{\be(\ta_* + 2\pi k)\}_{0\leq
k<N}$ has less than $N$ distinct values.  A twisted crossing point
$\ta_*$ is transverse when two requirments are met.  The first is met
when $\{\be(\ta_* + 2\pi k)\}_{0\leq k<N}$ has precisely $N - 1$
distinct elements.  Assuming now that the first requirement is met,
let $k \neq k'$ denote the two integers in $\{0,\ldots, N - 1\}$ for
which the value of $\be$ at $\ta_* + 2\pi k$ is the same as that at
$\ta_* + 2\pi k'$ agree.  The second requirement is then met when the
difference the locally defined function $\be(\ta + 2\pi k) - \be(\ta +
2\pi k')$ has non-zero derivative at $\ta =\ta_* $.

If the braid is such that its twisted crossing points are all
transverse, then those points where the $d\ta$ component of \eqref{eq26} has
less than N distinct values are in 1--1 correspondence with the set
of twisted crossing points.  Indeed, with this transversality
assumption, the final point in \eqref{eq16} guarantees that each point of the
one set is very close to precisely one point in the other.

Now, given the preceding comments, the final point in \eqref{eq16} has the
following implication: Let $\ta_*$ be a twisted crossing point and let
$k \neq k'\in \{0,\ldots, N - 1\}$ denote the unique pair for which
$\be(\ta_*+ 2\pi k) = \be(\ta_* + 2\pi k')$.  If
\begin{equation}
\al(\ta_*+ 2\pi k) - \al(\ta_* + 2\pi k')\quad \text{and}\quad      
\sin(\mu\ta_*+ 2\pi k/N) - \sin(\mu\ta_* + 2\pi k'/N)
\label{eq27}
\end{equation}
have the same sign, then such a $\ta_*$ corresponds to a unique
$\ta\in [0, 2\pi]$ near $\ta_*$ and a unique $r \in[1 - \de - \de^2, 1
- \de]$ where \eqref{eq27} has less than $N$ distinct points at $(r, \ta)$.
Furthermore, this then defines a $1-1$ correspondence between twisted
crossing points that obey \eqref{eq27} and points $(r, \ta)$ with $1 - \de -
\de^2 \leq r < 1 - \de$ where \eqref{eq27} has less than $N$ distinct values.

\NextPart6
This final part of the story explains how to use
information from the braid to compute the local sign of at the various
self intersection points of the Lagrangian.  For this purpose, agree
to orient the Lagrangian as a (multi-valued) graph over the $a_1$--plane,
where the latter is oriented by the form $r dr d\ta$.  This is to say that
the Lagrangian is to be viewed as the graph of the differential of
$f_{\ve}$ in \eqref{eq20}.

Now, suppose a transverse double point occurs in $L$ over a point in
the $a_1$--plane with coordinates $(r_*, \ta_*)$.  Thus, two sheets of
$L$ intersect at this point and so there exists a distinct pair $k,
k'\in \{0,\ldots, N- 1\}$ such that $df_{\ve}$ has the same value at
$(r_*, \ta_* + 2\pi k)$ and at $(r_*, \ta_* + 2\pi k')$.  This
understood, it then follows that the sign $(\pm 1)$ of this self
intersection is equal to minus the sign of the determinant of the
hessian at $(r_*, \ta_*)$ of the function
\begin{equation}
H \equiv f_{\ve}(r, \ta + 2\pi k) f_{\ve}(r, \ta + 2\pi k') .
\label{eq28}
\end{equation}
Apply this prescription to the $2(N - 1)$ self intersection points
described above in Part 3 to find that each has local intersection
number $-1$.  Indeed, to order $\calO(\ve)$, the function $H$ is the
same as $H_0 \equiv \chi_{\de} r^{\mu} (\sin(\mu\ta_* + 2\pi k) -
\sin(\mu\ta_* + 2\pi k'))$, and so small $\ve$ makes both the
differential and hessian of $H$ very close to those of $H_0$.  Thus,
small $\ve$ makes each of the relevant critical points of $H$ very
close to one of $H_0$ and it makes the signs of the corresponding
determinants agree if $H_0$'s determinant is not zero.  In this
regard, note that H0 has positive determinant at each relevant
critical point because each occurs where $\chi_{\de} r^{\mu}$ is
maximized.

Consider next the signs of the self intersection points that are
described above in Part 4.  In this regard, it follows from \eqref{eq24} that
when $\ve$ and $\de$ are very small, then the sign of the relevant
determinant is negative.  Indeed, this follows because the hessian in
question differs by $\calO(\ve\de)$ from a matrix having the form
$\ve\ch$, where $\ch$ is the symmetric matrix with zeros on the
diagonal and, in the notation from Part 4, with off diagonal entries
equal to $-w (\cos(\mu\ta_* + 2\pi k) - \cos(\mu\ta_* + 2\pi k'))$.
Thus, all of Part 4's self- intersection points have local
intersection sign equal to $+1$.

Turn at last to the self-intersection points that are described above
in Part 5.  Under the assumptions that all of the twisted crossing
points $\{\ta_*\}$ are non-degenerate and that $\ve$ and $\de$ are
both small, then the local intersection signs are determined as
follows: Suppose that a twisted crossing point $\ta_*$ determines a
self intersection point as described in Part 5.  Then, the local
intersection number for this intersection point is minus the product
of the sign of $\al(\ta_*+ 2\pi k)- \al(\ta_* + 2\pi k')$ with the
sign of the derivative at the point $\ta =\ta_* $ of $\be(\ta+ 2\pi k)
- \be(\ta + 2\pi k')$.  Here, $k$ and $k'$ are as given in \eqref{eq27}.  (The
latter all follows with the help of the final point in \eqref{eq16}.)

By the way, this sign can be interpreted as follows: View the triple
$(\al,\be,\ta)$ as the coordinates of a portion 
of the braid in $\bR^3$, and then
view $(\be,\ta)$ as the coordinates of the braid's 
projection into $\bR^2 \sub\bC$.  This done, 
then a twisted crossing point corresponds to a crossing
of strands as viewed via the direction defined by this projection.
Now, orient the the strand using the 1--form $d\ta$.  This understood, the
strand that corresponds near $\ta =\ta_*$  to the parameterization 
by $\ta\ra(\al(\ta + 2\pi k), \be(\ta + 2\pi k))$ 
passes on top of the other strand with respect to this projection 
when $\al(\ta_*+ 2\pi k) - \al(\ta_* + 2\pi k') > 0$ 
and passes under the other strand when 
$\al(\ta_*+ 2\pi k) - \al(\ta_* + 2\pi k') < 0$.  
This understood, the sign of the
corresponding self intersection point is positive when the crossing as
seen by this projection appears as in the following diagram:
\vspace{-.8cm}\begin{equation}
\hbox{
\setlength{\unitlength}{1mm}
\thicklines
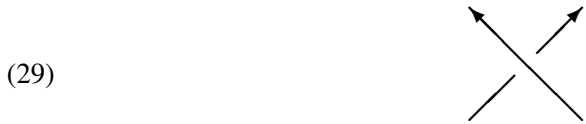
\begin{picture}(20,20)(0,5)
\put(20,0){\vector(-1,1){15}}
\put(5,0){\line(1,1){6}}
\put(14,9){\vector(1,1){6}}
\end{picture}
}\label{eq29}
\end{equation}\vspace{.4cm}
\end{ppart}

\section{Lagrangians and crossing changes}\label{sec3}

Suppose that two braids differ by a single strand crossing.  As
certain knot invariants can be characterized in terms of skein
relations, one might ask how the corresponding Lagrangians compare
with each other, and with that for the third braid in the skein
diagram.  To be more precise, suppose that the three braids are
identical except for their intersection with a fixed small ball in
$S^1\ti \bC$, and in this ball, the three braids correspond to the
following three pictures:\vspace{-.8cm}
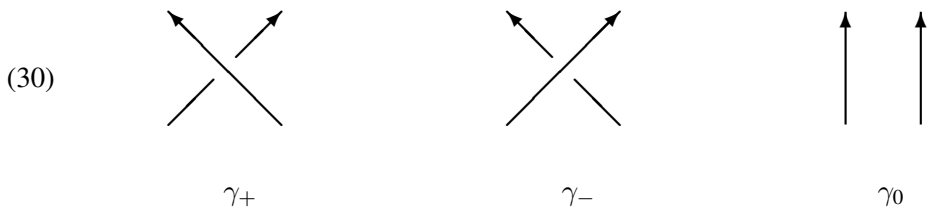
\begin{equation}
\lower1cm\hbox{
\setlength{\unitlength}{1mm}
\thicklines
\begin{picture}(80,30)(10,-5)
\put(20,0){\vector(-1,1){15}}
\put(5,0){\line(1,1){6}}
\put(14,9){\vector(1,1){6}}

\put(50,0){\vector(1,1){15}}
\put(65,0){\line(-1,1){6}}
\put(56,9){\vector(-1,1){6}}

\put(95,0){\vector(0,1){15}}
\put(105,0){\vector(0,1){15}}
\put(12,-10){$\ga_+$}
\put(57,-10){$\ga_-$}
\put(99,-10){$\ga_0$}
\end{picture}
}\label{eq30}
\end{equation}		

\bigskip\bigskip
This question is considered below when the Lagrangians for the braids
$\ga_{\pm}$ are such that they are given at large values of $r = |a_1|$ as the
graph of the differential of the appropriate version of the function
$f_*$ depicted in \eqref{eq19}.

Some conventions need setting to connect the pictures in \eqref{eq30} with
$f_*$.  For this purpose, suppose that when $\ga$ is one the braid in
one of the pictures in \eqref{eq30}, then the pair $(\al, \ta)$ give the $x$
and $y$ coordinates of the strands in \eqref{eq30}.  In this regard, the
convention is standard: The variable $x$ increases with horizontal
motion to the right in \eqref{eq30} and $y$ increases with vertical motion to
the top of the drawing in \eqref{eq30}.  This understood, then $\be$ should be
assumed to increase in the direction out of the paper but away from
the reader.  (Note that the drawing in \eqref{eq29} uses the different
convention where $(\be, \ta)$ are the coordinates of the projection.)

Now, to simplify notation, suppose that $\ta_* = 0$ is the value of
the $\ta$ coordinate where the projection in the $\ga_+$ diagram in
\eqref{eq30} has one strand pass over the other.  In this regard, note that
$\ga_+$ can be isotoped as a braid so that its parametrization at
values of $\ta$ near 0 is such that the under passing strand in the
$\ga_+$ diagram is described by $(\al_+(\ta) = \ta, \be_+(\ta)
\tau(\ta))$ where $\tau$ is a smooth, non-negative function of $\ta$
that is positive at $\ta = 0$ and vanishes where $\ta$ is near the top
and bottom of its implicit range in \eqref{eq30}.  To be precise, suppose that
this range for $\ta$ is $(-\ve, \ve)$ and that $\tau(\ta) = 0$ for
$|\ta| >\ve/2.$ Meanwhile, the over passing strand in the $\ga_+$
picture is described by $(\al_+(\ta + 2\pi k) = -\ta, \be_+(\ta + 2\pi
k) \equiv 0)$ where $k \in \{1,\ldots, N - 1\}$.

At the same time, $\ga_-$ can be isotoped as a braid so that for $\ta$
near 0, the parametrizing data $(\al_-, \be_-)$ has $(\al_-(\ta),
\be_-(\ta)) = (\ta, -\tau(\ta))$ and $(\al_-(\ta + 2\pi k), \be_-(\ta
+ 2\pi k))$ equal to $(-\ta, 0)$.  Thus, values of $\ta$ near zero in
both the $\ga_+$ and $\ga_-$ diagram parametrize the strand that
points up and to the right.

Having digested this notation, 
define the family $\{\ga_s: s\in[-1, 1]\}$ of maps from $S^1$ 
to $\bC$ as follows:  When $\ta\notin (-\ve, \ve)$, 
then $\ga_s(\ta) = \ga_+(\ta) = \ga_-(\ta)$.  
On the other hand, when $\ta\in (-\ve, \ve)$, then
\begin{equation}
(\al_s(\ta) = \ta, s \tau(\ta)).\label{eq31}
\end{equation}
This understood, then all positive s versions of $\ga_s$ define a
braid that is isotopic to $\ga_+$ while all negative s versions define
one that is isotopic to $\ga_-$.

Now consider the $\ga_s$ version, $L_s$, of the Lagrangian defined in \fullref{step8}
of \fullref{sec1} via the differential of \eqref{eq18}'s function $f_*$.  In
particular, the small, but positive s versions have a transversal
double point that is parametrized by $\ta = 0$ and $r = 2\tau(0)/s$ and
otherwise, no double points where $r \geq R_0$ with $R_0$ independent of $s$.
Meanwhile, the $s < 0$ versions of $L_s$ have no double points at all where
$r \geq R_0$.

Note next that even for $s = 0$, the definition given in \fullref{step8} of
\fullref{sec1} for $L_s$ makes perfectly good sense and describes a properly
immersed, Lagrangian disk in $\bC^2$.  In particular, $L_{s=0}$ can be assumed
to have solely transversal and isolated double points if the braid $\ga_+$
is chosen in a suitably generic fashion in the complement of the ball
pictured in \eqref{eq30}.  In any event, the non-zero versions of $L_s$ can be
assumed to converge in the $C^{\iy}$ topology as $s \ra 0$ to the $s = 0$
version.

In fact, these last conclusions about $L_s$ can be strengthened as
follows: The Lagrangian $L_0$ can be assumed to have the same double
points as all small and negative $s$ versions of $L_{-s}$ and be
isotopic by small Hamiltonian isotopies of $\bC^2$ to such $L_s$.  It
can also be assumed to intersect all spheres with radius greater than
$R_0$ transversely in a braid that is braid isotopic to $\ga_-$.
Meanwhile, Ls for small and positive s can be assumed to intersect all
spheres with radii between $R_0$ and $2\tau(0)/s$ transversely, and in
a braid that is also braid isotopic to $\ga_-$ even as it intersects
all spheres of radius greater than $2\tau(0)/s$ in a braid that is isotopic
to $\ga_+$.  Moreover, the portion of such a small and positive s version
of $L_s$ where $r$ is less than $2\tau(0)/s$ can be taken to be isotopic via
proper, Hamiltonian isotopy of $\bC^2$ to the same portions of the small,
but negative s versions of $L_s$.

This said about the $\ga_+$ and $\ga_-$ Lagrangians, what follows is a
description of a related Lagrangian, $L$, for the braid $\ga_0$ in
\eqref{eq30}.  For this purpose, fix some $r_0 \gg R_0$. Then $L$ has the
following properties:
\begin{equation}\mbox{}\label{eq32}\end{equation}
\begin{itemize}\vspace{-.5cm}
\item[$\bullet$]
The portion of $L$ where $r < r_0$ is isotopic via a Hamiltonian isotopy to
$L_0$, while the portion where $r > r_0$ is likewise isotopic to the $\ga_0$
version of the Lagrangian from $f_*$.  
\item[$\bullet$] 
In fact, given some positive $\ve > 0$, 
the $r < r_0 -\ve $ portion of $L$ can be taken equal to
$L_0$, while the $r > r_0 + \ve$ portion can be taken equal to the same
portion of the $\ga_0$ version of the Lagrangian from \fullref{step8} in \fullref{sec1}.
\item[$\bullet$]
$L$ is embedded where $r \geq R_0$, but this portion of $L$ is not a
cylinder, and thus not a multi-valued graph over the $a_1$--plane.
Rather, the $r \geq R_0$ portion of $L$ projects to the $a_1$ plane with
a single ramification point to account for the change in the topology
of its constant $r$ slices at $r = r_0$.
\end{itemize}

\medskip
The story on $L$ starts with a digression to provide a local model for
this ramification business. For this purpose, consider the locus in
$\bC^2$ where
\begin{equation}
a_1 - r_0 = a_2^2 .\label{eq33}
\end{equation}
Note that this locus defines a smooth, Lagrangian surface in $\bC^2$ whose
projection to the $a_1$--plane is $2-1$ save for the single critical point
that projects to $(r = r_0 , \ta = 0)$.
 	
To see how \eqref{eq33} models the desired behavior $\ta$ near zero, view the pair
$\al$ and $\be$ with $\ta$ as coordinates on 
$S^1 \ti \bC^2$.  Then, where $r$ is near $r_0$
and $\ta$ near 0, the equation in \eqref{eq33} has the schematic form
\begin{align}
\bullet&\quad 
\hbox{$r - r_0 = \al^2 - \be^2 + \cdots,$}\label{eq34}\\
\bullet&\quad 
\hbox{$ \ta = 2 r_0^{-1}\al\be + \cdots,$}\nonumber
\end{align}
where the `$\cdots$' signify terms that are $\calO(\al^4 + \be^4)$.  
	
This last equation understood, first fix $r$ at some value very near,
but less than $r_0$ and then view the resulting locus as a curve in
the $(\al, \ta, \be)$ version of $\bR^3$.  In particular, note that
view from the same vantage as that in \eqref{eq30} looks like the $\ga_-$
picture in \eqref{eq30}.  Meanwhile, the analogous view for the locus defined
by fixing $r$ near, but greater than $r_0$ in \eqref{eq34} looks like the
$\ga_0$ version of \eqref{eq30}.
	
Given the comments in the preceding paragraph, the task to construct $L$
as described in \eqref{eq32} is straightforward and left to the reader with
the hints to take $\al$ and $\be$ 
to be very small near $\ta = 0$ when comparing
with the description of the Lagrangian $L_0$.

\section{Lagrangians in $O(-1) \oplus O(-1)$}\label{sec4}

As in the introduction, let $O(-1) \ra \bP^1$ denote the degree $-1$,
holomorphic line bundle over the Riemann sphere.  The purpose of this
section is to describe how certain Lagrangians from \fullref{sec3} can be used to construct a 3--dimensional Lagrangian in the
K\"ahler manifold $O(-1) \op O(-1)$.  The construction starts with a 
2--dimensional Lagrangian, $L \sub \bC^2$, that is mapped to itself under
multiplication by $-1$ on $\bC^2$ and produces a 3--dimensional Lagrangian
in $O(-1) \op O(-1)$ that projects to the equator in $\bP^1$ with fiber $L$.

The symplectic form for the space $O(-1) \op O(-1)$ is a standard
K\"ahler form.  To view it, introduce the homogeneous complex
coordinates $(z, w) \in \bC^2 - \{0\}$ for $\bP^1$.  Thus, $(z, w)$
gives the same point in $\bP^1$ as $(\la z, \la w)$ when $\la\in \bC$
is not zero.  Now introduce the homogeneous coordinates $((z, w),
\eta_1, \eta_2)$ for $O(-1) \op O(-1)$ where now the latter and $((\la
z, \la w), \la^{-1} \eta_1, \la^{-1} \eta_2)$ give the same point.
This done, introduce the coordinates $\ze_i = (|z|^2 + |w|^2)^{1/2}
\eta_i$; the latter transform as $\ze_i \ra|\la| \la^{-1} \ze_i$ when
$(z, w) \ra (\la z, \la w)$.  In particular, the transformation for
each $\ze_i$ is unitary, so the norm $|\ze_i|$ gives a well defined
function.
	
Next, introduce 
\begin{equation}
\pi_i \equiv d\ze_i + A \ze_i\label{eq35}
\end{equation}
where $A$ is the connection 1--form 
\begin{equation}
A = (|z|^2 + |w|^2)^{-1}{\rm Im}\,(\bar{z} dz +\bar{w} dw); \label{eq36}
\end{equation}
thus $\pi_i$ transforms as does $\ze_i$ when $(z, w) \ra (\la z, \la
w)$ with $\la$ now any nowhere zero, complex valued function.  Letting
$u \ra (z = u, w = 1)$ denote the complex coordinate on the $w \neq 0$
portion of $\bP^1$, then $A =(|u|^2 + 1)^{-1}{\rm Im}\,(\bar{u} du)$
and $dA = (|u|^2 + 1)^{-2}{\rm Im}\,(d \bar{u}\wedge du)$.

With the $\{\pi_i\}$ in hand, the symplectic form 
on $O(-1) \op O(-1)$ is written
using the coordinate $u$ as
\begin{equation}
i 2^{-1}(\pi_1 \wedge\bar{\pi}_1 + \pi_2\wedge \bar{\pi}_2 
+ (|u|^2 + 1)^{-2} (a + |\ze_1|^2 + |\ze_2|^2) du \wedge d\bar{u});
\label{eq37}
\end{equation}
here $a > 0$ and $2\pi a$ gives the symplectic area of $\bP^1$.  For
reference in the subsequent discussion, note that this symplectic form
restricts to the $|u| = 1$ equator in $\bP^1$ as
\begin{equation}
i 2^{-1} \Si_i d(e^{i\varphi/2} \ze_i) 
\wedge d(e^{-i\varphi/2}  \bar{\ze}_i),\label{eq38}
\end{equation}
where $\va\in [0, 2\pi]$ is the argument of $u$.  

Now let $L \sub \bC^2$ denote an immersed, Lagrangian surface that is mapped
to itself by the action of multiplication by $-1$ on $\bC^2$.  With \eqref{eq38}
understood, it follows directly that
\begin{equation}
M \equiv \{(u = e^{i\va}, \ze_1, \ze_2): 
(e^{i\va/2} \ze_1, e^{i\va/2} \ze_2) \in L\}\label{eq39}
\end{equation}
is an immersed Lagrangian in $O(-1) \op O(-1)$.  As an abstract
manifold, $M$ is diffeomorphic to the quotient of $S^1 \ti L$ by the
action of $\bZ/2$ that sends $(\va, z_1, z_2)$ to the point $(\va +
\pi, -z_1, -z_2)$; a diffeomorphism here is provided by the map that
sends the equivalence class of $(\va, z_1, z_2)$ to $(\va, \ze_1 =
e^{-i\va/2} z_1, \ze_2 = e^{-i\va/2} z_2)$.  Note that if $L$ is
embedded in $\bC^2$, then $M$ is embedded in $O(-1) \op O(-1)$ and if
$L$ is immersed with transverse double points, then $M$ is immersed
too.  However, the immersion of the latter is not transverse since it
is a union of circles.

By the way, this construction is identical to that given by Equation
(5.3) of \cite{LMV} when applied to the hyperk\"ahler rotation of the
zero locus of a holomorphic function in $\bC^2$.

Additonal examples come from the construction above in \fullref{sec1}.  In
particular, suppose that the braid $\ga$ has the following property: There
exists $k \in \{0,\ldots, N - 1\}$ such that
\begin{equation}
\ga(\ta + 2\pi (k + 1/2)) = - \ga(\ta)\label{eq40}
\end{equation}
at each $\ta\in [0, 2\pi]$.  This condition asserts that the braid is
mapped to itself by the action of multiplication by $-1$ on $\bC^2$.
For example, \eqref{eq40} holds when $N$ is odd and $\ga = e^{i\ta/N}$.

In any event, if $N$ is odd and if \eqref{eq40} holds, then the constructions
in \fullref{sec1} produce Lagrangians from the braid $\ga$ that are mapped
to themselves by the $-1$ action on $\bC^2$.  For example, the connect
sum of any knot with itself can be represented by a braid with this
property.

\section{The view from $T^*S^3$}\label{sec5}

As remarked at the outset, Gopakumar and Vafa came to their conjecture
by applying `t Hooft's ideas to certain string theories on $T^*S^3$.
This application suggested a duality between these string theories on
$T^*S^3$ and others on $O(-1) \op O(-1)$.  As this duality has an explicit
geometric basis, the Lagrangians just constructed in $O(-1) \op O(-1)$ can
be viewed from the perspective of $T^*S^3$.  Such is the purpose of this
final section.

The geometric basis for afore-mentioned duality is simply that both
$O(-1) \op O(-1)$ and $T^*S^3$ can be viewed as resolutions of the singularity
at the origin in $\bC^4$ of the zero locus of a certain quadratic
polynomial.  To be more precise, introduce complex coordinates 
$(c_1,\ldots, c_4)$ for $\bC^4$.  This done, the polynomial in question is
\begin{equation}
\wp = c_1c_2 - c_3c_4 .\label{eq41}
\end{equation}
The total space of $O(-1) \op O(-1)$ then maps onto $\wp^{-1}(0)$ via a
holomorphic map that is one to one off of the zero section and
collapses the latter to the origin in $\bC^4$.  The map in question sends
the homogeneous coordinates $(z, w,\eta_1, \eta_2)$ to
\begin{equation}
(c_1 = z \eta_2, c_2 = w \eta_1, c_3 = z \eta_1, c_4 = w \eta_2) .\label{eq42}
\end{equation}
Meanwhile, $T^*S^3$ maps to $\wp^{-1}(0)$ as follows:  
First, take two copies of $\bR^4$ and use  $y \equiv (y_1,\ldots, y_4)$ 
to denote a point in the first and 
$v \equiv(v_1,\ldots, v_4)$ for a point in the second.  
This done, identify the complement of the zero section in $T^*S^3$ with
the subset of  $\bR^4 \ti \bR^4$ where $|y| = |v| \neq 0$ 
and $\Si_k y_kv_k = 0$.  Here, the convention taken is that the  
assignment of $|y|^{-1} y \in S^3$ to $(y, v)$ 
defines the projection to $S^3$.
Now, identify $\bR^4 \ti \bR^4$  with $\bC^4$ via 
\begin{align}
\bullet&\quad\hbox{
 $c_1 = y_1 + i v_1 - i (y_2 + i v_2)$,}\label{eq43}\\ 
\bullet&\quad\hbox{
 $c_2 = y_1 + i v_1 + i (y_2 + i v_2)$,}\nonumber\\
\bullet&\quad\hbox{
 $c_3 = -(y_3 + i v_3) + i (y_4 + i v_4)$,}\nonumber\\
\bullet&\quad\hbox{
 $c_4 = y_3 + i v_3 + i (y_4 + i v_4)$.}\nonumber 
\end{align}
This map sends the complement of the zero section of $T^*S^3$
diffeomorphically onto the complement of the origin in $\wp^{-1}(0)$,
and it extends in the obvious way as a smooth map from $T^*S^3$ onto
the whole of $\wp^{-1}(0)$ that sends the zero section to the origin
in $\bC^4$.

Now, the preceding describes the `correspondence diagram'
\begin{equation}
O(-1) \op O(-1) \wp^{-1}(0) \ra T^*S^3 , \label{eq44}
\end{equation}
where both arrows are diffeomorphisms from the complement of the
corresponding zero sets to the complement of the origin.

This correspondence gives the following convoluted map from $(S^1 \ti
\bC^2)/\{\pm 1\}$ to $\bC^2$: First, embed this space 
in $O(-1) \op O(-1)$ by the map that sends a point labeled by a unit
length complex coordinate $u$ for $S^1$ and pair $(z_1, z_2)$ of
complex coordinates for $\bC^2$ to the point with the homogeneous
coordinates
\begin{equation}
(z = u,w = 1,\eta_1=2^{-1/2} u^{-1/2} z_1, 
\eta_2 = 2^{-1/2} u^{-1/2} z_2) .\label{eq45}
\end{equation}
Clearly, the image of this map fibers over the equator in $\bC\bP^1$
with fiber $\bC^2$, and the intersection of the image with the zero
section of $O(-1) \ti O(-1)$ is the image of $S^1 \ti \{0\}$.  Next,
use the left arrow in \eqref{eq44} (thus, \eqref{eq42}) to identify the complement of
the $S^1 \ti \{0\}$ in $(S^1 \ti \bC^2)/\{\pm 1\}$ with a subset of
the complement of the origin in $\wp^{-1}(0) \sub \bC^4$.  This done,
use the inverse of the right arrow in \eqref{eq44} (the inverse of \eqref{eq43}) to
identify the complement of $S^1 \ti \{0\}$ in $(S^1 \ti \bC^2)/\{\pm
1\}$ with a subset of the complement of the zero section in
$T^*S^3$. Finally, project the latter to $\bC^2$ using the projection
from $T^*S^3$ to $\bC^2$ that sends $(y, v)$ to the point with the
coordinates $(y_1 + i y_2, y_3 + i y_4)$.  This map extends as a
smooth map from $(S^1 \ti \bC^2)/\{\pm 1\}$ to $\bR^4$ sending $S^1
\ti \{0\}$ to the origin.  Amusingly, the map just described is very
simple when written with the hyperk\"ahler rotated coordinates $(a_1,
a_2)$ in \eqref{eq2}.  Indeed, this map sends $(u, (a_1, a_2))$ to
\begin{equation}
2^{-1} (u^{-1/2} a_1, u^{-1/2} a_2) .\label{eq46}
\end{equation}
The preceding has the following implications:Let $L \sub \bC^2$ be a
Lagrangian surface that is mapped to itself via multiplication on
$\bC^2$ by $-1$.  Construct from L the 3--dimensional Lagrangian $M =
(S^1 \ti L)/\{\pm 1\}$ in $O(-1) \op O(-1)$ as described in \eqref{eq39}.
This done, use the correspondences in \eqref{eq44} to identify the complement
of $M$'s intersection with the zero section with a subset, $M^*$, in
$T^*S^3$.  Finally, map $M^*$ to $\bC^2$ via the map from $T^*S^3$
that assigns $(y_1 + i y_2,y_3 + i y_4)$ to $(y, v)$.

The result is a smooth map from $(S^1 \ti L)/\{\pm 1\}$ to $\bC^2$
that simply rotates $L$ as in \eqref{eq46}.  In particular, if $L$ intersects
some 3--sphere about the origin as a knot, then for each fixed $u \in
S^1$, the corresponding image of $(u \ti L)$ in $\bC^2$ intersects the
concentric half radius sphere as a rotated image of the same knot.

\bibliographystyle{gtart}
\bibliography{link}

\end{document}